\documentclass[fleqn]{mat01}
\usepackage{times,mathtimy,amssymb,latexsym,amsbsy}
\begin{document}

\setcounter{page}{67}
\firstpage{67}

\newtheorem{theore}{Theorem}
\renewcommand\thetheore{\arabic{section}.\arabic{theore}}
\newtheorem{theor}[theore]{\bf Theorem}
\newtheorem{definit}[theore]{\rm DEFINITION}
\newtheorem{lem}[theore]{Lemma}
\newtheorem{pot}[theore]{Proof of Theorem}

\newtheorem{theoree}{Theorem}
\renewcommand\thetheoree{{\it \roman{theoree}}}
\newtheorem{case}[theoree]{\it Case}

\newtheorem{cla}{\it Claim}

\def\egztheorem{\trivlist\item[\hskip\labelsep{{\bf EGZ Theorem.}}]}
\def\theorr{\trivlist\item[\hskip\labelsep{{\bf Theorem.}}]}
\def\cauchy{\trivlist\item[\hskip\labelsep{{\it Cauchy--Davenport inequality}}]}
\def\openpro{\trivlist\item[\hskip\labelsep{{\it Open Problem.}}]}
\def\claim{\trivlist\item[\hskip\labelsep{{\it Claim.}}]}
\def\conjecture{\trivlist\item[\hskip\labelsep{{\it Conjecture}}]}

\def\zn{{\Bbb Z}_n}
\def\zp{{\Bbb Z}_p}

\title{On the structure of $\pmb{p}$-zero-sum free sequences and its
application to a variant of Erd\"{o}s--Ginzburg--Ziv theorem}

\markboth{W~D~Gao, A~Panigrahi and R~Thangadurai}{On the structure of $p$-zero-sum free sequences}

\author{W~D~GAO, A~PANIGRAHI$^{*}$ and R~THANGADURAI$^{*}$}

\address{Department of Computer Science and Technology, University of
Petroleum, Changping Shuiku Road, Beijing 102200, China\\
\noindent $^{*}$School of Mathematics, Harish-Chandra Research Institute,
Chhatnag Road, Jhusi, Allahabad~211~019, India\\
\noindent E-mail: wdgao\_1963@yahoo.com.cn; anupama@mri.ernet.in; thanga@mri.ernet.in}

\volume{115}

\mon{February}

\parts{1}

\pubyear{2005}

\Date{MS received 8 September 2004}

\begin{abstract}
Let $p$ be any odd prime number. Let $k$ be any positive integer such
that $2\leq k\leq \left[\frac{p+1}3\right]+1$. Let $S =
(a_1,a_2,\ldots,a_{2p-k})$ be any sequence in ${\Bbb Z}_p$ such that
there is no subsequence of length $p$ of $S$ whose sum is zero in $\zp$.
Then we prove that we can arrange the sequence $S$ as follows:
\begin{equation*}
S = (\underbrace{a, a, \ldots, a}_{u \ {\rm times }}, \underbrace{b, b,
\ldots, b}_{v \ {\rm times}}, a_1', a_2', \ldots, a_{2p-k-u-v}')
\end{equation*}
where $u\geq v$, $u+v\geq 2p-2k+2$ and $a-b$ generates $\zp$. This
extends a result in \cite{gao10} to all primes $p$ and $k$ satisfying
$(p+1)/4+3\leq k\leq (p+1)/3+1$. Also, we prove that if $g$ denotes the
number of distinct residue classes modulo $p$ appearing in the sequence
$S$ in $\zp$ of length $2p-k$ $(2\leq k\leq [(p+1)/4]+1)$, and $g\geq
2\sqrt{2}\sqrt{k-2}$, then there exists a subsequence of $S$ of length
$p$ whose sum is zero in $\zp$.
\end{abstract}


\maketitle

\section{Introduction}

Let $n$ be any positive integer. Let $S = (a_1,a_2,\ldots,a_\ell)$ be a
sequence (possibly with repetition) in the cyclic group of order $n$
(denoted by ${\Bbb Z}_n$) of length $\ell$. We call a subsequence $T =
(b_1,b_2,\ldots,b_r)$ of $S$ to be zero-sum subsequence if
$b_1+b_2+\cdots+ b_r = 0$ in $\zn$.

In 1961, Erd\"{o}s--Ginzburg--Ziv proved the following theorem (which we
call the EGZ theorem).

\begin{egztheorem}\hskip -.5pc {\bf \cite{egz}.}\ \ {\it Given a sequence $S$ in
$\zn$ of length $2n-1${\rm ,} one can extract a zero-sum subsequence of
length $n$ in $\zn$.}\vspace{.7pc}
\end{egztheorem}

The EGZ theorem is tight in the following sense. If
\begin{equation*}S = (\underbrace{0,0,
\ldots,0}_{n-1 \ {\rm times }}, \underbrace{1,1,\ldots,1}_{n-1 \ {\rm
times}})
\end{equation*}
is a sequence in $\zn$ of length $2n-2$, then $S$ does not have a
zero-sum subsequence of length $n$.

Many authors studied the characterization of the above extremal example.
In particular, Yuster and Peterson \cite{yp1} and independently
Bialostocki and Dierker \cite{bd} proved that any sequence $S$ in $\zn$
of length $2n-2$ having no zero-sum subsequence of length $n$ will be of
the form
\begin{equation*}
S = (\underbrace{a,a,\ldots,a}_{n-1 \ {\rm times }}, \underbrace{b,b,
\ldots,b}_{n-1 \ {\rm times }}),
\end{equation*}
where $a\ne b \in \zn$.

Also, Flores and Ordaz \cite{fo1} proved the following result of this
nature. Suppose $S$ is any sequence in $\zn$ of length $2n-3$ such that
$S$ has no zero-sum subsequence of length $n$. Then there exists $a, b
\in \zn$ such that $\zn$ is generated by $b-a$ and $a$ appearing $n-1$
times in $S$ and one of the following conditions hold: (i) $b$
appearing exactly $n-2$ times; (ii) $b$ appearing exactly $n-3$ times in
$S$ and also, $2b -a$ appearing exactly once in $S$.

In 1996, Gao \cite{gao10} proved the generalization of the above two
results as follows.\vspace{.7pc}

\begin{theorr}\hskip -.5pc {\bf \cite{gao10}.}\ \ {\it Let $n$ be any
positive integer. Let $k$ be any positive integer such that $2\leq k\leq
\left[\frac{n+1}4\right]+2$. Let $S = (a_1,a_2,\ldots,a_{2n-k})$ be any
sequence in ${\Bbb Z}_n$ such that there is no subsequence of length $n$
of $S$ whose sum is zero in $\zn$. Then we can re-arrange the sequence
$S$ as follows{\rm :}
\begin{equation*}
S = (\underbrace{a, a,\ldots,a}_{u \ \mbox{ times }}, \underbrace{b,
b,\ldots,b}_{v \ {\rm times}}, a_1',a_2',\ldots,a_{2n-k-u-v}')
\end{equation*}
where $u\geq v, u+v\geq 2n-2k+2$ and $a-b$ generates $\zn$.}\vspace{.7pc}
\end{theorr}

One of our main theorems in this article is to extend the above result
to all primes $p$ and integer $k$ for the range $\frac{p+1}4+3\leq k
\leq \frac{p+1}3+1$. This extension is meaningful for all large primes
$p$. Also, we shall study the problem of how many distinct residue
classes modulo $p$ occur in those sequences of length $2p-k$ in $\zp$
having a zero-sum subsequence of length $p$ in it. Before we state our
main theorems, we shall fix up notations as\break follows.

For every integer $1\leq k\leq \ell$, define
\begin{equation*}
\sum_{k}(S)=\left\{a_{i_1}+a_{i_2}+\cdots+a_{i_k} | 1\leq i_1 < i_2 <
\cdots < i_k \leq \ell\right\}
\end{equation*}
and $\sum(S) = \cup_{k=1}^\ell \sum_k(S)$. For any subsequence $T =
(b_1, b_2,\ldots,b_r)$ of $S$, we let $\sigma(T) = \sum_{i=1}^rb_i$. We
denote $ST^{-1}$ by the deleted sequence $R$ which is obtained from $S$
by deleting the elements of $T$. Also, if $S =
({\underbrace{a,a,\ldots,a}_{r \ {\rm times}}}, b_1,b_2,\ldots)$,
then we write $S = (a^r,b_1,b_2,\ldots)$. For any $b\in {\Bbb Z}_n$, we
denote by $b+S$ the sequence $(b+a_1,b+a_2,\ldots, b+a_\ell)$. For every
$x \in {\Bbb Z}_n$, define $\overline{x}$ to be the least positive
inverse image under the natural homomorphism from the additive group of
integers ${\Bbb Z}$ onto ${\Bbb Z}_n$. For example, $\overline{0} = n$.
If $A\subset {\Bbb Z}_n$, then we denote the cardinality of $A$ by
$|A|$. If $A$ is a sequence in ${\Bbb Z}_n$, we denote the length of $A$
by $|A|$ (same notation as the cardinality). For any $g\in{\Bbb Z}_n$,
we define $v_g(S)$ by the number of times $g$ appears in $S$. Also, we
define $h = h(S) = \max_{g\in{\Bbb Z}_n}v_g(S)$. Gao \cite{gao10}
introduced the following definition.

\pagebreak

\begin{definit}$\left.\right.$\vspace{.5pc}

\noindent {\rm Let $S = (a_1,a_2,\ldots,a_\ell)$ and $T = (b_1,
b_2,\ldots,b_\ell)$ be two sequences in ${\Bbb Z}_n$ of length $\ell$.
We say that $S$ is equivalent to $T$ (written as $S\sim T$) if there
exist an integer $c$ coprime to $n$, an element $x\in {\Bbb Z}_n$, and
a permutation $\pi$ of $\{1,2,\ldots,\ell\}$ such that $a_i =
c(b_{\pi(i)} -x)$ for every $i = 1,2, \ldots, n$. Clearly, $\sim$ is
an equivalence relation; and if $S\sim T$, then $0\in\sum_n(S)$ if and
only if $0\in \sum_n(T)$.}
\end{definit}

In this article, we shall prove theorems~3.1 and 3.2.

\setcounter{section}{3}
\setcounter{theore}{0}
\begin{theor}[\!] Let $p$ be any odd prime number. Let $k$ be any
positive integer such that $2\leq k\leq \left[\frac{p+1}3\right]+1$. Let
$S = (a_1,a_2,\ldots,a_{2p-k})$ be any sequence in ${\Bbb Z}_p$ such
that $0\not\in \sum_p(S)$. Then
\begin{equation*}
S \sim (0^u, 1^v, a_1',a_2',\ldots,a_{2p-k-u-v}'),
\end{equation*}
where $u\geq v$ and $u+v\geq 2p-2k+2$.
\end{theor}

Using the information in Theorem~3.1, we consider the following problem of
variant of EGZ theorem as follows. Before we state our theorem, we
recall the following definition which was introduced in \cite{bl} and
state the known results.

\setcounter{section}{1}
\setcounter{theore}{1}
\begin{definit}$\left.\right.$\vspace{.5pc}

\noindent {\rm Let $n,k$ be positive integers, $1\leq k\leq n$.
Denote by $f(n,k)$ the least positive integer $g$ for which the following
holds: If $S = (a_1, a_2, \ldots, a_g)$ is a sequence of elements of
${\Bbb Z}_n$, the cyclic group of order $n$, of length $g$ such that the
number of distinct $a_i$'s is equal to $k$, then there are $n$ indices $i_1,
i_2, \ldots, i_n$ belonging to $\{1, 2, \ldots, g\}$ such that $a_{i_1}+
a_{i_2}+ \cdots+ a_{i_n}\break = 0$.}
\end{definit}

\vspace{.2pc}
\begin{theorr}{\it We have
\begin{enumerate}
\renewcommand\labelenumi{\rm \arabic{enumi}.}
\leftskip -.15pc
\item $f(n,k) \leq 2n-1$ for all $n$ and for all $1\leq k\leq n$ {\rm
(}By EGZ theorem{\rm )}.\vspace{.3pc}

\item $f(n,n) =\left\{ \begin{array}{ll} n, &{\rm if } \ n \ {\rm
is~odd}\\ n+1, &{\rm if } \ n \ {\rm is~even } \end{array}\right.$
{\rm \cite{ggp}}.\vspace{.3pc}

\item $f(n,k) = n+2${\rm ,} for all $n\geq 5$ and $1+n/2 < k \leq n-1$
{\rm \cite{br,ggp}}.\vspace{.3pc}

\item $f\left(n, \frac{n}2+1 \right) = n+3$ for all $n\in 2{\Bbb N}$
{\rm \cite{gghhlp}}.\vspace{.3pc}

\item $f(n,k) = 2n -((k-1)/2)^2-1$ for all $n\geq (k-1)^2-4$ for an odd
$k\geq 5$ {\rm \cite{wan}}.\vspace{.3pc}

\item $f(n,k) = 2n - k(k-2)/4 -1$ for all $n \geq k(k-2)-4$ for an even
$k\geq 6$ {\rm \cite{wan}}.\vspace{.3pc}

\item $f(n,2) = 2n-1, f(n,3) = 2n-2$ and $f(n,4) = 2n -3$ for all $n$
{\rm \cite{bl}}.\vspace{.3pc}

\item $f(n,k) \leq 2n -k+1$ for all $2\leq k\leq n$ {\rm \cite{ham}}.\vspace{.3pc}

\item $f(p,k) \leq 2p -3k+11$ for all $5\leq k \leq (p+15)/3$
{\rm \cite{at}}.
\end{enumerate}}
\end{theorr}

Other than these results many authors (for instance \cite{gg},
\cite{bl} and \cite{biol}) consider some lower bounds for $f(n,k)$ for
various $k$.

In this article, we shall prove the following result.

\setcounter{section}{3}
\setcounter{theore}{1}
\begin{theor}[\!] Let $p$ be any odd prime number. Let $k$ be any
positive integer such that $2\leq k\leq \left[\frac{p+1}3\right]+1$.
Then $f(p, \ell) \leq 2p-k$ for all $\ell \geq 2\sqrt{2} \sqrt{k-2}$.
\end{theor}

\setcounter{section}{1}
\section{Preliminaries}

We shall start this section with a well-known fundamental inequality of
subsets as follows.

\begin{cauchy}\hskip -.3pc {\rm \cite{cau,dav}.}\ \
Let $p$ be any prime number. Let $A_1, A_2, \ldots, A_t$ be
non-empty subsets of ${\Bbb Z}_p$. Then\vspace{-.2pc}
\begin{equation*}
|A_1+A_2+\cdots+A_t| \geq \min\left\{p, \sum_{i=1}^t|A_i| -t+1\right\}.
\end{equation*}

$\left.\right.$\vspace{-.8pc}

\end{cauchy}

\setcounter{section}{2}
\setcounter{theore}{0}
\begin{theor}[\!]\hskip -.5pc {\bf \cite{ben}.}\ \
Let $n$ and $k$ be any positive integers such that $n-2k\geq 1$. If $S =
(a_1, a_2, \ldots, a_{n-k})$ is a sequence in ${\Bbb Z}_n$ such that
$0\not\in\sum(S)${\rm ,} then there exists $a\ne 0\in {\Bbb Z}_n$ which
appear at least $n-2k+1$ times in $S$.
\end{theor}

The following Theorem is crucial for the proof of Theorem~3.1.

\begin{theor}[\!] Let $p$ be any prime number and $1\leq k \leq p-2$.
Let $S$ be a sequence in ${\Bbb Z}_p$ of length $p+k$. If $0\not\in
\sum_p(S)${\rm ,} then $h(S)\geq k+1$.
\end{theor}

\begin{proof}
When $k=1$, the result follows from the Pigeon hole principle. So, we
can assume that $k\geq 2$. If possible, we assume that $h(S) \leq k$.
Then, we can distribute the elements of $S$ into a union $A_{1} \sqcup A_{2}
\sqcup \cdots \sqcup A_{k}$, so that in each $A_{i}$, an element occurs only
once. By the Cauchy--Davenport theorem, we see that
\begin{align*}
\left|\sum_{i=1}^{k}A_i\right| &\geq \min\left\{p, \sum_{i=1}^{k}|A_i| -
k+1\right\}\\[.3pc]
&= \min\{p, p+k-k+1 = p+1\} = p.
\end{align*}
Therefore, $A_1+A_2+\cdots+A_{k} = \zp$. In particular, $\sigma(S)\in
\sum_{k}(S)$. Without loss of generality we shall assume that $\sigma(S)
= a_1+a_2+\cdots+a_{k}$. Then we have $a_{k+1}+a_{k+2}+\cdots+ a_{p+k} =
0$ which implies $0\in \sum_p(S)$ as $|S| = p+k$. This contradicts the
assumption that $0\not\in \sum_p(S)$. Therefore, $h(S)\geq k+1$.
$\hfill\Box$
\end{proof}

\begin{theor}[\!]\hskip -.5pc {\bf \cite{gao11}.}\ \ Let $n$ be any positive
integer. Let $1\leq k\leq \left[\frac{n+1}3\right]${\rm ,} and let $S$
be a sequence in ${\Bbb Z}_n$ of length $n-k$ such that
$0\not\in\sum(S)$. Then
\begin{equation*}
S \sim (1^{n-2k+1}, x_1,x_2,\ldots,x_{k-1})
\end{equation*}
with $\sum_{i=1}^{k-1}\overline{x_i} \leq 2k-2$.
\end{theor}

\begin{lem} Let $p$ be any odd prime and $1\leq k\leq \left[ \frac{p+1}3
\right]$. Let $S =(1^{p-2k+1}, x_1, x_2, \ldots${\rm ,} $x_{k-1})$ be a sequence
in ${\Bbb Z}_p\backslash\{0\}$ of length $p-k$ such that
$\sum_{i=1}^{k-1} \ \overline{x_i} \ \leq 2k-2$. Then{\rm ,}\, for any
$x\in {\Bbb Z}_p$ satisfying $p-2k+1\leq \ \overline{x} \ \leq p-2k+1+
\sum_{i=1}^{k-1}\ \overline{x_i}${\rm ,} there exists a subsequence $T$
of $S$ such that $|T| \geq p-2k+1$ with $\sigma(T) = x$.
\end{lem}

\begin{proof}
Let $x\in{\Bbb Z}_p$ such that $p-2k+1\leq \ \overline{x} \ \leq
p-2k+1+\sum_{i=1}^{k-1}\ \overline{x_i}$. If $x = p-2k+1$, then $x =
\sum_{i=1}^{p-2k+1} 1$ and we are done; otherwise, we have
\begin{equation*}
p-2k+2\leq \ \overline{x} \ \leq p-2k+1+\sum_{i=1}^{k-1} \
\overline{x_i} \ \leq p-1.
\end{equation*}
Therefore, we have $1\leq \ \overline{x} - (p-2k+1) \ \leq
\sum_{i=1}^{k-1} \ \overline{x_i}$.

\begin{claim} For any positive integer $k$, if $S' = (x_1,x_2,\ldots,
x_{k})$ be a sequence in ${\Bbb Z}_p\backslash\{0\}$ such that $|S'| =
k$ and $\sum_{i=1}^k\overline{x_i}\leq 2k$, then, for every $x\in {\Bbb
Z}_p$ satisfying $1\leq \overline{x} \leq \sum_{i=1}^k\overline{x_i}$,
either $x\in\sum(S')$ or $x+1\in \sum(S')$.

If the claim is proven, then, we get, either $x -(p-2k+1)$ or
$x-(p-2k+1)+1$ in $\sum((x_1,x_2,\ldots,x_{k-1}))$. That is, either $x =
\underbrace{1+1+ \cdots+1}_{p-2k+1}+y$ or $x =
\underbrace{1+1+\cdots+1}_{p-2k}+y$ where $y\in
\sum((x_1,x_2,\ldots,x_{k-1}))$. So, to end the proof of this lemma, it
is enough to prove this claim.\vspace{.7pc}
\end{claim}

When $k = 1,2$, the claim is trivially true. So, we let $k\geq 3$.
Assume the result is true for $k-1$ and we shall prove for $k$. If
necessary by renaming the indices, without loss of generality, we can
assume that $S' = (x_1,x_2,\ldots,x_k)$ with $\overline{x_1}\leq
\overline{x_2}\leq \cdots \leq \overline{x_k}$. Suppose
$\overline{x_{k-1}} = 1$. Then, we have $x_1 = x_2=\cdots=x_{k-1}=1$. As
$\sum_{i=1}^k\overline{x_i} \leq 2k$, we see that $\overline{x_k} \leq
2k - (k-1) = k+1$. Therefore, we see that
\begin{equation*}
\sum(S') = \left\{\begin{array}{ll}
\{1,2,\ldots,x_k+k-1\}, &{\rm if } \ \overline{x_k} \leq k,\\[.4pc]
\{1,2,\ldots,x_k+k-1\}\setminus\{k\},  &{\rm if } \ \overline{x_k} =
k+1 \end{array}\right.
\end{equation*}
which clearly implies the claim. Thus, now, we can assume that $2\leq \
\overline{x_{k-1}} \ \leq \ \overline{x_k}$. If $\overline{x} \ \leq \
\overline{x_k} \ + \sum_{i=1}^{k-2}\ \overline{x_i}$, then by induction,
either $x$ or $x+1$ in $\sum((x_1,x_2, \ldots,x_{k-2},x_k))$, and we are
through; otherwise, we have, $\overline{x_k} \ + \ \sum_{i=1}^{k-2}\
\overline{x_i} \ \leq \ \overline{x} \ < \ \sum_{i=1}^k \
\overline{x_i}$. Therefore, we have
\begin{equation*}
k-2 \ \leq \overline{x_k}-\overline{x_{k-1}}+ \sum_{i=1}^{k-2} \
\overline{x_i} \ \leq \ \overline{x-x_{k-1}} \ \leq \ \overline{x_k} \ +
\ \sum_{i=1}^{k-2} \ \overline{x_i}.
\end{equation*}
Therefore, by the induction hypothesis, we see that either $x-x_{k-1}$
or $x-x_{k-1}+1$ in $\sum((x_1, x_2, \ldots,x_{k-2},x_k))$ and hence, we
have either $x$ or $x+1$ in $\sum(S').\hfill\Box$
\end{proof}

\section{Proof of Theorems 3.1 and 3.2}

\setcounter{theore}{0}
\begin{pot}{\rm
Let $S$ be a sequence in ${\Bbb Z}_p$ of length $2p-k$ where $2\leq
k\leq [\frac{p+1}3]+1$. Given that $0\not\in\sum_p(S)$. Without loss of
generality we can assume that $0$ (if necessary, by translating by an
element) appears maximum number of, say $u$, times in $S$. By
Theorem~2.2, it is clear that $u\geq p-k+1$. Therefore, $S =
(0^{u},a_1,a_2,\ldots,a_{2p-k-u})$ where $a_i\in {\Bbb Z}_p\backslash
\{0\}$. Let $S_1 = (a_1,a_2,\ldots,a_{2p-k-u})$ be a subsequence of $S$.
Since $u\geq p-k+1$, we have $2p - k -u\leq 2p -k-p+k-1 = p-1$. That is,
$|S_1|\leq p-1$. Let $|S_1| = p-m$ for some positive integer $m$. Note
that $p-m+u= 2p-k$ which implies $u+k-p = m$. As $0\not\in\sum_p(S)$, we
have $u\leq p-1$. Therefore, $m = u+k-p\leq p-1+k-p$\break $= k-1$.

If $0\not\in \sum(S_1)$, then by Theorem 2.1, we know that there exists
an element $a\in {\Bbb Z}_p\backslash\{0\}$ such that $v_a(S_1)\geq
p-2m+1$. Therefore, $S = (0^u,a^v, b_1,b_2,\ldots,b_{2p-k-u-v})$ and
$2p-k-u-v\leq m-1\leq k-2$ which implies $2p-2k+2\leq u+v$ and we\break are
done.

Thus, we can assume that $0\in\sum(S_1)$. Let $W$ be the maximal
zero-sum subsequence of $S_1$ of length $w$. Moreover, since $0\not\in
\sum_p(S)$ and $S_1$ is a sequence in $\zp\backslash\{0\}$ and $u\geq
p-k+1$, we have
\begin{equation}
2\leq w \leq  p-u-1 \ \Longrightarrow \ 2\leq w\leq k-2.
\end{equation}
Also note that $k+u+w \geq k+p-k+1+w\geq p+1$. Put $\ell = k+u+w -p$.
Therefore, $2p -k-u-w = p -\ell$. By the definition of $W$, we have
$0\not\in\sum(S_1W^{-1})$ and $|S_1W^{-1}| = p-\ell$. Let $T =
S_1W^{-1}$. Also, by the inequality (1), we see that $\ell = k+u+w-p
\leq k+u+p-u-1-p = k-1\leq \left[\frac{p+1}{3}\right]$. Therefore, by
Theorem 2.3, we see that
\begin{equation*}
T \sim (1^{p-2\ell+1},x_1,x_2,\ldots,x_{\ell-1}) \quad \hbox{ and }
\quad \sum_{i=1}^{\ell-1}\overline{x_i}\leq 2\ell -2.
\end{equation*}
Thus, the given sequence $S = 0^uS_1 = 0^uTW$ is equivalent to the
following sequence:
\begin{equation*}
S \sim (0^{u},1^{p-2\ell+1},x_1,x_2,\ldots,x_{\ell-1}, z_1, z_2, \ldots,
z_w)
\end{equation*}
where all the $x_i\ne 0$ satisfying
$\sum_{i=1}^{\ell-1}\overline{x_i}\leq 2\ell -2$ and $W \sim
(z_1,z_2,\ldots,z_w)$ is the maximal zero-sum subsequence of $S_1$.

Without loss of generality, we shall replace `$\sim$' by `$=$' above. Also,
we denote the number of $1$'s appearing in the sequences
$(x_1,x_2,\ldots, x_{\ell-1})$ and $(z_1, z_2, \ldots, z_w)$ by $r$ and
$t$ respectively. Put $v = p-2\ell +1+r+t$.

To end the proof of this theorem, it is enough to prove that $u+v\geq
2p-2k+2$.

If $2\leq \ \overline{z_i} \ \leq p-2\ell+1$ for some $i$ satisfying
$1\leq i \leq w$, then as there are $p-2\ell+1$ number of $1$'s in $T$,
we can write $z_i = \sigma(L_1)$ where $L_1 = (1^{\overline{z_i}})$ with
$|L_1|\geq 2$. If $p-2\ell+2\leq \overline{z_i}\leq p-2\ell+1+
\sum_{j=1}^{\ell-1}\ \overline{x_i} $ holds for some $1\leq i\leq w$,
then by Lemma 2.4, there exists a subsequence $L_1$ of $T$ such that
$z_i = \sigma(L_1)$ and $|L_1| \geq 2$. By letting $W_1 = L_1Wz_i^{-1}$,
we see that $\sigma(W_1) = 0$ and $|W_1| \geq w+1$ which contradicts the
maximality of $W$. Hence
\begin{equation}
p-1\geq \ \overline{z_i} \ \geq p-2\ell + 2+\sum_{i=1}^{\ell-1} \
\overline{x_i} \quad \hbox{ for each } \ z_i\ne 1.
\end{equation}
Since $\sum_{i=1}^{\ell-1} \ \overline{x_i} \ \leq 2\ell-2$, we have
\begin{equation}
2\ell-2\leq \sum_{i=1}^{\ell-1} \ \overline{x_i} +r.
\end{equation}
Therefore, by the inequalities (2) and (3), we get
\begin{equation}
p-1\geq \ \overline{z_i} \ \geq p-r \quad \hbox{ for each } \ z_i\ne 1.
\end{equation}
By rearranging the indices and renaming them, if necessary, we can assume
that for $0\leq q \leq w$, we have
\begin{equation}
\overline{z_i} \ \ne 1 \quad \hbox{ for } \ 1\leq i\leq q \quad \hbox{
and } \quad \overline{z_i} = 1 \quad \mbox{ for } \ q+1\leq i\leq w.
\end{equation}
\pagebreak

\begin{case}
$(w = 2)$
\end{case}

\noindent In this case, by the definition of $W$, we have $z_1+z_2 = 0$.
Therefore, there are two cases, namely, $z_1 = 1$ and $z_2 = -1$ or
$z_1\ne 1$ and hence $z_2\ne -1$. When $z_1\ne 1$, by the inequality
(4), we have $p-2\geq \overline{z_1}\geq p-r$ and in particular, we have
$r\geq 2$. Since $2\leq r \leq p-2\ell+1$, we have a zero-sum
subsequence $Z = (\overline{z_1}, 1^{p-\overline{z_1}})$ which has
length $ > 2$ which is a contradiction to the maximality of $W$. Thus,
$z_1 = 1$ and $z_2 = -1$. In this case, $v\geq 2p-k-u-\ell+r$.
Therefore, $u+v = u+2p-k-u-\ell+r \geq 2p -k-(k-1)+r = 2p -2k+1+r\geq
2p-2k+3$. We are done in this case.

\begin{case}{\rm
$(w\geq 3)$}
\end{case}

\noindent Since $W$ is a zero-sum sequence, $q\ne 0$. So, we have $1\leq
q\leq w$. When $q =1$, from the inequality (4), we get
\begin{equation}
p-1 \geq  \ \overline{z_1} \  \geq p-r.
\end{equation}

When $q = 2$, we have
\begin{align*}
2p-2 &\geq \ \overline{z_1} \ + \ \overline{z_2} \ \geq 2p -2r\geq 2p -
2(\ell-1)\geq 2p-2(k-2)\\[.3pc]
&= 2p -2k+4.
\end{align*}
Since $k\leq \left[\frac{p+1}3\right]+1$, it is clear that $p\geq 3k-4$
and hence
\begin{equation*}
2p-2\geq \ \overline{z_1} \ + \ \overline{z_2} \ > p \ \Longrightarrow
p-2\geq \ \overline{z_1} \ + \ \overline{z_2} \ - p > p- 2\ell+2.
\end{equation*}
Therefore, it follows that
\begin{equation*}
p-2\geq \ \overline{z_1+z_2} \ > p-2\ell+2.
\end{equation*}
If $\ \overline{z_1+z_2} \ \leq p-2\ell+1+\sum_{i=1}^{\ell-1} \
\overline{x_i}$, then $z_1+z_2 = \sigma(L_2)$ for some subsequence $L_2$
of $T$ with $|L_2|\geq p-2\ell+1$ (by Lemma 2.4). If we let $W_2 =
L_2Wz_1^{-1} z_2^{-1}$, then $\sigma(W_2)= 0$ and $|W_2|=|L_2|+w-2 \geq
p-2\ell+1+w -2 = w+p-(2\ell+1) > w$ (as $\ell < k-1\leq (p+1)/3$) which
contradicts the maximality of $W$. Therefore, we have
$\overline{z_1+z_2} \geq \sum_{i=1}^{\ell-1} \ \overline{x_i} \ +
p-2\ell+2$. Thus, by the inequality (3), we have
\begin{equation}
p-2\geq \ \overline{z_1+z_2} \ \geq p-r.
\end{equation}

Now, we shall assume that $q\geq 3$. Let $a = \min\{q,w-2\}$. Then we
claim the following.

\begin{cla}{\rm For $q\geq 3$ and for every $s = 1,2,\ldots,a$, we have
\begin{equation*}
p -s \geq \left(\sum_{i=1}^s \ \overline{z_i} \ \right) - (s-1)p =
\overline{\sum_{i=1}^sz_i} \geq p-r.
\end{equation*}}
\end{cla}

By inequalities (4) and (7), the Claim 1 is true for $s = 1$ and $2$.
Now, we shall assume that claim 1 is true for $s-1$ and we prove for
$s$. By the inequality (4) and induction hypothesis, we have
\begin{align*}
p-s &\geq \left(\sum_{i=1}^s \ \overline{z_i} \ \right) - (s-1)p =
\left(\sum_{i=1}^{s-1} \ \overline{z_i} \ \right) - (s-2)p + \
\overline{z_s} \ -p\\[.3pc]
&\geq p-r-r= p-2r\geq p-2\ell+2 \geq p-2k+4 > 0.
\end{align*}
Hence,
\begin{equation*}
p-2\ell+2\leq  \ \overline{\sum_{i=1}^sz_i}.
\end{equation*}
If $ \ \overline{\sum_{i=1}^sz_i} \ \leq p-2\ell+1+\sum_{i=1}^{\ell-1} \
\overline{x_i}$, then by Lemma 2.4, there exists a subsequence $L_3$ of
$T$ with $|L_3|\geq p-2\ell+1$ such that $\sum_{i=1}^sz_i =
\sigma(L_3)$. If we let $W_3 = L_3Wz_1^{-1}z_2^{-1}\ldots z_s^{-1}$,
then we get $\sigma(W_3) = 0$. Since $w\leq k-2$, $\ell\leq k-1$ and
$p\geq 3k-1$, we have
\begin{align*}
|W_3| &= w+|L_3|-s\geq w + p-2\ell+1 -(w-2)\\[.3pc]
&\geq w + p - 2k + 4 + 1 - k + 4 = w + p - (3k - 9) > w.
\end{align*}
This contradicts the fact that $W$ is the maximal zero-sum subsequence
of $S_1$. Therefore, we have
\begin{equation*}
p-s\geq \overline{\sum_{i=1}^sz_i} \ \geq p-2\ell+2+\sum_{i=1}^{\ell-1}
\ \overline{x_i}
\end{equation*}
and by the inequality (3), we get Claim 1.

\begin{cla}{\rm $q \leq w -2$.}
\end{cla}

\noindent Assume, on the contrary that $q\geq w-1$. Then $q = w-1$ or $q
= w$. If $q = w-1$, then we have $p -(w-2)\geq
\overline{z_1+z_2+\cdots+z_{w-2}} \ \geq p-r$, $p-1\geq
\overline{z_{w-1}} \ \geq p-r$ and $z_w = 1$. Therefore,
\begin{align*}
2p &> 2p -w +2 \geq \ \overline{z_1+z_2+\cdots+z_{w-2}} \ + \
\overline{z_{w-1}} \ + \ \overline{z_w}\\[.3pc]
&\geq 2p-2r+1 > p
\end{align*}
which is a contradiction to $\sigma(W)=0$. Hence $q\ne w-1$.

If $q = w$, then $p -(w-2)\geq \ \overline{z_1+z_2+\cdots+z_{w-2}} \
\geq p-r$, \ \ $p-1\geq \ \overline{z_{w-1}}, \ \ \overline{z_w} \ \geq
p-r$. Therefore,
\begin{align*}
3p &> 3p -w \geq \ \overline{z_1+z_2+\cdots+z_{w-2}} \ + \
\overline{z_{w-1}} \ + \ \overline{z_w}\\[.3pc]
&\geq 3p -3r \geq 3p -3(k-2)\geq 3p -3k+6 = 2p +p-3k+6> 2p,
\end{align*}
(as $r\leq \ell-1\leq k-2$ and $p\geq 3k-4$) which is also a
contradiction to $\sigma(W) =0$. Hence $q\ne w$. Thus Claim 2 is
true.

From Claims 1 and 2, we see that $s$ varies from $1$ to $q$. Since we
have $p -s \geq p-r$ which implies $r\geq s$. In particular, when $s = q$,
we get
\begin{equation}
q\leq r.
\end{equation}
But by the definition of $q$, we have $q = w-t$ which implies that $w =
q+t$. Therefore, by the inequality (8), we have $r+t \geq q+t= w$. Thus
\begin{align*}
u+v &= u+p -2\ell+1+r+t\geq u+p-2\ell+1+w\\[.3pc]
&= 2p -k-(\ell-1)\geq 2p -2k+2,
\end{align*}
as desired.} $\hfill\Box$
\end{pot}

\begin{pot}{\rm Let $S$ be a given sequence in $\zp$ of length $2p-k$.
Suppose the number of distinct residue classes appearing in $S$ is
$g\geq 2\sqrt{2}\sqrt{k-2}$. If possible, we assume that $0\not\in
\sum_p(S)$. Then by Theorem~3.1, $S = 0^u T W$ (notations as in the
proof of Theorem~3.1). Now, we shall count the number of distinct
residue classes modulo $p$ appearing in $T$ and in $W$ separately.

We recall that $T = (1^{p-2\ell+1}, x_1,x_2,\ldots, x_{\ell-1})$ with
$\sum_{i=1}^{\ell-1} \overline{x_i}\leq 2\ell-2$ and $r =
v_1((x_1,x_2,\ldots,x_{\ell-1}))$. Also, $W = (z_1, z_2, \ldots, z_q,
\underbrace{1, 1, \ldots, 1}_{w-q \ \mbox{ times }})$ where $z_i\ne 1$.
Note that by Claim~2 of Theorem~3.1, we have $1\leq q\leq w-2$ and by
(6) and (8) we have $q\leq r$ and $r\geq 2$.

Let $g_1$ (respectively, $g_2$) denote the number of distinct residue
classes modulo $p$ appearing in $T$ (respectively, in $W$). Thus,
including $0$, the total number of distinct residue classes modulo $p$
appearing in $S$ is $g=g_1+g_2+1-1 = g_1+g_2$ because the residue $1$ is
calculated twice in $g_1$ and $g_2$. So, to end the proof of this
theorem, it is enough to estimate $g=g_1+g_2$.

Since $\sum_{i=1}^{\ell-1}\overline{x_i}\leq 2\ell-2$ and $r$ number of
$1$'s appearing in $(x_i)$s, we have
\begin{align*}
&1+2+\cdots+g_1 \leq 2\ell -2-(r-1)\\[.3pc]
&g_1^2+g_1 \leq 4\ell -4-2(r-1)\leq 4\ell-4-2 = 2(\ell-3).
\end{align*}
Therefore, since $\ell \leq k-1$, we have
\begin{equation}
g_1^2+g_1 \leq 2(k-4) \ \Longrightarrow \ g_1 \leq \sqrt{2} \sqrt{k-4}<
\sqrt{2} \sqrt{k-2}.
\end{equation}
Now, note that\, $\overline{-z_i} = p -\overline{z_i}$. Therefore by
Claim~1 of Theorem~3.1, we get $\sum_{i=1}^q\overline{-z_i}\leq r$. Thus,
\begin{equation*}
1+2+\cdots+g_2\leq r \Longrightarrow g_2 \leq \sqrt{2r}.
\end{equation*}
Since $r\leq \ell-1\leq k-2$, we have
\begin{equation}
g_2\leq \sqrt{2k-4} = \sqrt{2}\sqrt{k-2}.
\end{equation}
Thus, from the inequalities (9) and (10) and counting $0$, we have
\begin{equation*}
g_1+g_2 < \sqrt{2}\sqrt{k-2}+\sqrt{2}\sqrt{k-2} = 2\sqrt{2}\sqrt{k-2},
\end{equation*}
a contradiction. Hence the theorem.} $\hfill\Box$
\end{pot}

We shall end this section with the following open problems.\vspace{.7pc}

\begin{openpro} Let $n$ and $k$ be two positive integers such that
$k\leq n-2$. Determine the constant defined by
\begin{equation*}
h(n,k)=\min\{h(S) \ | \  |S|=n+k\},
\end{equation*}
where $S$ runs over all sequences in $\zn$ of length $n+k$ such that
$0\not\in \sum_n(S)$.\vspace{.7pc}
\end{openpro}\pagebreak

It is proved in \cite{bd} and \cite{yp1} that $h(n,n-2)=n-1$ and proved
in \cite{fo1} that $h(n,n-3)=n-1$. Theorem~2.2 shows that $h(p, k)
\geq k+1$ for all $1 \leq k \leq p-2$. The main result in \cite{gao10}
implies that $h(n,k)\geq k+1$ whenever $n-[(n+1)/4]-1\leq k\leq n-2$. It
is natural to ask if $h(n,k)\geq k+1$ for every positive integer $n$ and
every $k$ such that $1\leq k\leq n-2$. However, the answer is `no' in
general. Recently, in \cite{gaothaw} we provided a counter example for
$k$ satisfying $p \leq k\leq n/p-2$. We conjectured the following.\vspace{.7pc}

\begin{conjecture}\hskip -.3pc \cite{gaothaw}.\ \
Let $n > 1$ be any positive integer and let $p$ be the smallest prime
divisor of $n$. Let $k$ be an integer such that $k\geq (n/p)-1$.
Then $h(n,k) \geq k+1$.\vspace{.7pc}
\end{conjecture}

In \cite{gaothaw}, it is proved that Conjecture 1 is true for $n =
p^\ell$ for any prime $p$. Also, it is not known whether Conjecture 1 is
true for $k < p/3$.

\section*{Acknowledgements}

The first author is supported by NSFC with grant Nos 19971058 and
10271080.

\end{document}